\documentclass[10pt,a4paper,twoside]{article}

\usepackage{amsfonts,amssymb,amsmath,amsthm,graphics}

\usepackage{fancyhdr}
\usepackage{titlesec}
\titleformat{\section}
  {\normalfont\large\bfseries}{\thesection}{0.5em}{}
\newcommand{\addperiod}[1]{#1.}
\titleformat{\subsection}[runin]
  {\normalfont\bfseries}{\thesubsection}{0.25em}{\addperiod}

\renewcommand\thesection{\arabic{section}.}
\renewcommand\thesubsection{\thesection\arabic{subsection}}

\theoremstyle{plain}
\newtheorem*{definition}{\bf Definition}
\newtheorem*{remark}{\it Remark}
\newtheorem*{remarks}{\it Remarks}
\newtheorem{theorem}{\bf Theorem}
\newtheorem{proposition}{\bf Proposition}

\newtheorem*{example}{\it Example}

\newcommand\mathringring[1]{
  {\mathop{\kern0pt #1}\limits^{\vbox to-1.85ex{
   \kern-2ex \hbox to 0pt{\hss\normalfont\kern.1em \r{}\kern-.45em
  \r{}\hss} \vss }}}}

\newcommand{\der}[2]{\frac{\partial#1}{\partial#2}}
\newcommand{\dder}[3]{\frac{\partial^2#1}{\partial#2 \partial#3} }
\newcommand{\cC}{\mathcal{C}}
\newcommand{\cD}{\mathcal{D}}
\newcommand{\cM}{\mathcal{M}}
\newcommand{\cQ}{\mathcal{Q}}
\newcommand{\cR}{\mathcal{R}}
\newcommand{\bR}[1]{\mathbb{R}^{#1}}

\begin{document}

\title{\Large \bf 
Conservation of `moving' energy
\\
in nonholonomic systems with affine constraints
\\
and integrability of spheres 
on rotating surfaces\footnote{
This work is part of the research projects {\it
Symmetries and integrability of nonholonomic
mechanical systems} of the University of Padova and PRIN {\it Teorie
geometriche e analitiche dei sistemi Hamiltoniani in dimensioni finite
e infinite}.}
}

\author{\sc Francesco Fass\`o\footnote{\footnotesize Universit\`a di Padova,
Dipartimento di Matematica, Via Trieste 63, 35121 Padova, Italy.
Email: {\tt fasso@math.unipd.it} } 
\  and
Nicola Sansonetto\footnote{\footnotesize Universit\`a di Padova,
Dipartimento di Matematica, Via Trieste 63, 35121 Padova, Italy.
Email: {\tt sanson@math.unipd.it} }\ \footnote
{Supported by the Research Project
{\it Symmetries and integrability of nonholonomic
mechanical systems} of the University of Padova.}
}

\date{}
\maketitle
\centerline{\small (\today)}

{\small 
\begin{abstract}
\noindent
Energy is in general not conserved for mechanical nonholonomic systems with
affine constraints. In this article we point out that, nevertheless,
in certain cases, there is a  modification of the energy that is
conserved. Such a function coincides with the energy of the system
relative to a different reference frame, in which the constraint is
linear.
After giving sufficient conditions for this
to happen, we point out the role of symmetry in this mechanism.
Lastly, we apply these ideas to prove that the motions of
a heavy homogeneous solid sphere that rolls inside a convex 
surface of revolution in uniform rotation about its vertical figure
axis, are (at least for certain parameter values and in open regions
of the phase space) quasi-periodic on tori of dimension up to three.

\vskip3mm
{\scriptsize
\noindent
{\bf Keywords:} Nonholonomic mechanical systems, Conservation of
energy, Rolling rigid bodies, Symmetries and momentum maps,
Integrability.
\vskip1mm
\noindent
{\bf MSC:} 70F25, 37J60, 37J15, 70E18
}

\end{abstract}
}

\section{Introduction}

This paper is ultimately addressed to the class of mechanical systems
formed by a rigid body subject to the nonholonomic constraint of
rolling without sliding on a surface which moves in a preassigned
way.  This type of nonholonomic constraints are affine, not linear,
in the velocities. Consequently, even if the system is
time-independent (which may easily be the case in presence of
symmetries of body and surface, and if the latter moves at uniform
speed) energy need not be conserved. In fact, while energy is
conserved in time-independent nonholonomic systems with constraints
that are linear in the velocities (see e.g. \cite{pars,NF}),
the same is not always true if the constraints are affine in the velocities
(see e.g. \cite{marle2003,kobayashi-oliva} and Section 3.5 below). A
simple example of this situation is the classical system formed by a
sphere that rolls without sliding on a table that rotates uniformly,
studied by Pars \cite{pars}, Neimark and Fufaev~\cite{NF} and others,
in which the energy is not conserved.

The starting point of this paper is the observation that for a
nonholonomic system with affine constraints,  even if the energy is
not conserved, there might exist a modification of the energy---that
may be interpreted as the energy of the system relative to a different
reference frame and for this reason will be called {\it moving
energy}---that is conserved. The reason underneath this fact
is simply that, in a moving reference frame in which the surface is
at rest, the constraint is linear. Therefore, if the system happens
to be time-independent relatively to such a moving frame, its
energy relative to that frame is conserved. And the condition of
time-independence is easily verified in presence of symmetries.

We study the existence of conserved moving energies in Section 3. For
simplicity, instead of changing the reference frame with respect to
which the system is described, we pass to a moving system of
coordinates. After illustrating the mechanism on the well known 
example of a sphere on a turntable, we give sufficient conditions for
the existence of a conserved moving energy (Theorem \ref{thm1}). Even
though these conditions might appear very special, we relate them to
the existence of symmetries: we assume that a group acts in
configuration space and show that certain invariance properties of
the system lead to a conserved moving energy (Theorem~\ref{thm2}).
Interestingly, this conserved function is the sum of two functions
that, at variance from what would happen in a holonomic system, are
not conserved: the energy of the system and a component of the
momentum map of the lifted action. There is here a connection with the
failure of Noether theorem in nonholonomic mechanics, that we discuss.

Lastly, in Section 4 we apply the results of Section 3 to the study
of the system formed by a heavy homogeneous solid sphere that rolls
without sliding inside a convex surface of revolution, which rotates
around its (vertical) figure axis with constant speed~$\Omega$. 
When $\Omega=0$ this system is known
to be integrable, with dynamics quasi-periodic on tori of dimension
up to three \cite{hermans,zenkov}. If $\Omega\not=0$ the system
is $\mathrm{SO(3)}\times
S^1$-invariant and, so far, it was only known that its
four-dimensional reduced system admits two first integrals and a
conserved measure, and thus that it is integrable by quadratures
\cite{BMK2002}. By exploiting the existence of an additional first
integral given by a moving energy, much stronger
integrability results can be obtained. Here we give a first, general
result in this direction. Specifically, using essentially a
continuity argument from the case $\Omega=0$, we prove that,
for small $\Omega$, there is an open nonempty subset of the reduced
phase space in which the reduced dynamics is periodic, and
correspondingly an open nonempty subset of the unreduced phase space in
which the unreduced dynamics is quasi-periodic on tori
of dimension up to three (Theorem~\ref{thm3}).

Even though our primary interest is toward time-independent
nonholonomic systems, the need of considering time-dependent
coordinates forces us to work in the time-dependent context. This
somewhat complicates the notation. In the hope of keeping the
complexity to a minimum, we adopt a Lagrangian description on the
extended phase space  of time-dependent mechanical nonholonomic systems,
which is quickly described in Section 2.

For general introductions to nonholonomic mechanics see e.g. \cite{pars,NF,
cortes, bloch,  marle2003,benenti,CDS}; the
time-dependent case is treated, using the formalism of jet bundles, in
\cite{pagani91,sarlet95,sarlet96,pagani97}.  Throughout the
paper all manifolds and maps are assumed to be smooth, and all vector
fields are assumed to be complete. 

\vfill\eject

\section{Nonholonomic systems with affine constraints}

\subsection{Time-independent nonholonomic systems with affine
constraints}

First we briefly recall the time-independent case, mainly
to fix some notation. The
starting point is a Lagrangian system with $n$-dimensional
configuration manifold $Q$ and Lagrangian \hbox{$L:TQ\to\bR{}$}, that we
assume to be regular; hence, in each set of local bundle coordinates
$(q,\dot q)$ the matrix $\dder L {\dot q}{\dot q} (q)$ is everywhere
invertible. 

An affine nonholonomic constraint consists in the prescription that
the kinematic states of the system belong to the fibers $(\cM_0)_q$,
$q\in Q$, of an affine distribution $\cM_0$ on $Q$, that we
assume to have constant rank $r>1$ and to be nonintegrable. Thus,
there  are a vector field $\xi_0$ on $Q$ and a nonintegrable
distribution $\cD_0$ of constant rank $r$ on $Q$ such that
$$
  (\cM_0)_q = \{ v_q\in T_qQ : \, 
                 v_q-\xi_0(q)\in (\cD_0)_q \} 
  \quad \forall q \in Q \,.
$$
Of course, given $\cM_0$ and $\cD_0$, $\xi_0$ is defined up to a section of
$\cD_0$. $\cD_0$ will be said to be the distribution associated to
$\cM_0$. The affine distribution $\cM_0$ may also be regarded as a
submanifold $M_0$ of $TQ$, which is in fact an affine subbundle of
$TQ$. We call $M_0$ the {\it constraint submanifold}. The case of
linear constraints is recovered when $\xi_0=0$; the constraint
manifold is thus a linear subbundle of $TQ$.

We assume that the nonholonomic constraint is `ideal',
that is, that it satisfies d'Alembert principle (see e.g.
\cite{AKN-3d-ed,pagani91, marle2003}): at each $q\in Q$, the set of
reaction forces that the constraint can exert is (an appropriate jet
extension of) the annihilator $(\cD_0)_q^\circ \subset T^*_qQ$ of
the  fiber $(\cD_0)_q\subset T_qQ$ of the distribution
$\cD_0$ associated to the constraint submanifold.
It is well known that, under this assumption, there is a
unique choice of the reaction force as a function
$$
   R_{L,M_0}:M_0\to \cD_0^\circ 
$$
such that the restriction to $M_0$ of Lagrange equations with
the reaction force defines a vector field on $M_0$ (see e.g.
\cite{ago, benenti}). We denote this vector field on
$M_0$ as $X_{L,Q,M_0}$ and call it {\it time-independent nonholonomic
system}, {\it with affine constraints} if $M_0$ is an affine subbundle of
$TQ$ and {\it with linear constraint} if $M_0$ is a linear subbundle of
$TQ$.

\subsection{Definitions and notation for the time-dependent case} 
In order to consider time-dependent nonholonomic systems
we pass to the extended phase space. In doing so, we
need a number of definitions that we collect in this section.
Let $Q$ be an $n$-dimensional manifold.

An {\it$m$-dimensional extended submanifold $M$ of $TQ$} is an
$(m+1)$-dimensional submanifold of the extended phase space $TQ\times
\bR{}$ of the form
\begin{equation}\label{Mext}
  M = \{ M_t\times \{t\} \,:\; t\in\bR{} \}
  \,.
\end{equation}
Thus, for each $t$, $M_t$ is a submanifold of $TQ$ of dimension $m$.
The reason for the use of the term `extended', instead of the perhaps
more natural `time-dependent', is that we need to treat both
time-dependent  and time-independent cases
within the same context. We say that the extended submanifold $M$
is {\it time-independent} if
$$
  M=M_0\times \bR{}
$$
for a given submanifold $M_0$ of $TQ$, or equivalently if  $M_t=M_0$
for all $t$, and that it is {\it time-dependent} otherwise. 

If all the $M_t$'s in (\ref{Mext}) are linear subbundles of $TQ$, then the
extended submanifold $M$ is an {\it extended linear subbundle} of $TQ$. If
they are all affine subbundles, then $M$ is an {\it extended affine
subbundle} of $TQ$. Obviously, we regard extended subbundles as
special cases of extended affine subbundles.

{\it An extended distribution $\cD$ on $Q$} is a distribution on $Q\times
\bR{}$ with fibers
$$
  \cD_{(q,t)} = (\cD_t)_q \oplus \{0\}
$$
where, for each $t$, the $(\cD_t)_q$ are the fibers of a
distribution $\cD_t$ on $Q$.  If all the distributions $\cD_t$ have
rank $r$, then we say that $\cD$ has rank $r$. We say
that $\cD$ is nonintegrable if (some at least of) the
distributions $\cD_t$ are nonintegrable. An extended
distribution $\cD$ on $Q$ of rank~$r$ generates an
extended linear subbundle $D=\{D_t\times \{t\} : t\in \bR{}\}$ of
$TQ$ of dimension $n+r$, with  $D_t=\{v_q :  q\in Q, v_q \in
(\cD_t)_q\}$, and vice versa. $\cD$ is {\it time-independent} if $\cD_t=\cD_0$ for all
$t$.

An {\it extended vector field} on $Q$ is a vector field $\xi$ on $Q\times
\bR{}$ whose $\bR{}$-component is identically equal to $0$, namely
$$
  \xi(q,t)=\xi_t(q)+0\,\partial_t 
$$
with each $\xi_t$ a vector field on $Q$. $\xi$ is {\it
time-independent} if
$\xi_t=\xi_0$ for all $t$.

If $\cD$ is an extended distribution on $Q$ of rank $r$ and $\xi$ is an
extended vector field on~$Q$, then
\begin{equation*}
  \cM = \cD + \xi 
\end{equation*}
is an {\it affine extended distribution} of rank $r$ on $Q$.
Thus $\cM$ has fibers $\cM_{(q,t)} = (\cD_t)_q+\xi_t(q)$, or
$\cM_{(q,t)} = (\cM_t)_q\oplus \{0\}$ with $\cM_{t} = \cD_t +\xi_t$.
$\cD$ is called the {\it extended distribution associated to $\cM$}.
$\cM$ can be regarded in an obvious way as an extended affine subbundle
$M$ of $TQ$ of dimension $n+r$.  $\cM$ is {\it time-independent} if so are $\xi$
and $\cD$.

Finally, a {\it dynamical system} on an extended submanifold $M$
of $TQ$ is a vector field on $M$ whose $\bR{}$-component is
identically equal to $1$, that is 
$$
   X(v_q,t) = X_t(v_q)+\partial_t
$$
with each $X_t$ a vector field on $M_t$. (The difference with respect
an extended vector field is that now time does not stay constant, which is
necessary for the dynamics). $X$ is {time-independent} if $X_t=X_0$ for all $t$.

\subsection{Time-dependent nonholonomic systems}
We start now from a Lagrangian system with $n$-dimensional
configuration manifold $Q$ and time-dependent regular Lagrangian
$L:TQ\times\bR{}\to\bR{}$. The time dependency of the Lagrangian
accounts, for instance, for the presence of time-dependent holonomic
constraints.\footnote{A large class of time-dependent holonomic constraints for
systems of $N$ material points can be modelled in this way. After the
choice of a reference frame, that provides a (time-dependent) identification of
physical 3-space with $\bR3$, a time-dependent
holonomic constraint is given by a time-dependent embedding of a manifold $Q$ into the
configuration space $(\bR3)^N$ of the unconstrained system, and the
Lagrangian is the restriction of the Lagrangian of the unconstrained
system to the resulting time-dependent, extended submanifold.} 

We add now the nonholonomic constraint that its kinematical states
belong to an extended affine subbundle
$M=\{M_t\times\{t\}:t\in\bR{}\}$ of $TQ$ of dimension $n+r$, for some
$1<r<n$, that we call the {\it extended constraint submanifold}. This
extended affine subbundle corresponds to an extended affine distribution $\cM = \cD +
\xi$ of rank $r$ on $Q$. We assume that the associated extended
distribution $\cD$ is nonintegrable.

The condition of `ideality' of the constraint now means that, at each $t$
and $q$, the set of reaction forces that the constraint can exert is
(a jet extension of) the annihilator $(\cD_t)_q^\circ \subset T^*_qQ$ of the 
fiber $(\cD_t)_q$ \cite{pagani91, marle2003}, and implies that there
is a unique choice of the reaction force as a function
$$
   R_{L,M}:M\to \cD^\circ 
$$
such that the restriction to $M$ of Lagrange equations
with the reaction force defines a dynamical system 
$X_{L,Q,M}$ on $M$. Here, $\cD^\circ$ is the extended codistribution on
$Q\times\bR{}$ with fibers $(\cD_t)_q^\circ\oplus\{0\}$. 

\begin{definition}
Let $L:TQ\times \bR{}\to\bR{}$ a regular Lagrangian and $M$ an
extended affine subbundle of $TQ$. 
\begin{list}{}
{\leftmargin2em\labelwidth1.2em\labelsep.5em\itemindent0em
\topsep0.5ex\itemsep-0.2ex}
\item[i.] The dynamical system $X_{L,Q,M}$ on
$M$ is called the {\rm nonholonomic system with affine constraints}
(or, shortly, the {\rm nonholonomic system}) with Lagrangian $L$ and
extended constraint manifold $M$. 
\item[ii.] If $M$ is an extended linear subbundle of $TQ$ then we say that $X_{L,Q,M}$ has
{\it linear constraints.} 
\item[iii.] 
If $L$ and $M$ are time-independent, then we say that $X_{L,Q,M}$ is
{\rm time-independent}.
\end{list}
\end{definition}

In the time-independent case we will routinely identify $X_{L,Q,M_0}$ and
$X_{L,Q,M_0\times\bR{}}$, and, depending on the context, we will
regard the Lagrangian $L$ as defined either on $TQ$ or on
$TQ\times\bR{}$.

In (possibly time-dependent) bundle coordinates $(q,\dot q)$ in $TQ$,
the fibers of the distributions $\cD_t$ on $Q$ are the kernels of a
$(q,t)$-dependent $k\times n$ matrix $S(q,t)$ that has everywhere
rank $k$, with $k=n-r$: 
$$
   (\cD_t)_q =\{\dot q\in T_qQ \,:\; S(q,t)\dot q=0\} \,.
$$
Thus, $\dot q\in (\cM_t)_q$ if and only if $\dot q= \xi_t(q)+u$ with
$u\in\ker S(q,t)$, that is, if and only if  $S(q,t)[\dot
q-\xi_t(q)]=0$. It follows that, for each $t$ and $q$,
the affine subspace $(\cM_t)_q$ of $T_qQ$ is described by
\begin{equation}
\label{Sands=0}
   S(q,t)\dot q + s(q,t) = 0
\end{equation}
with
$
   s(q,t) = -S(q,t) \xi_t(q) \in \bR{k} 
$.
Of course, only $\ker S$ is uniquely defined, not $S$, $\xi$ and~$s$.
In coordinates, the  annihilator of the fiber $(\cD_t)_q$ is the range
of the matrix $S(q,t)^T$, and the reaction force $R_{L,M}(q,\dot
q,t)\in \mathrm{range}[S(q,t)^T]$. 

\subsection{Time-dependent diffeomorphisms and conjugation of
nonholonomic systems}

In order to implement time-dependent coordinate changes, we need to
consider (lifted) diffeomorphisms of the configuration space that
depend on time, and use them to transform nonholonomic systems. (We
may use now the expression `time-dependent', instead of `extended', because 
we will never need to consider
`time-independent' time-dependent change of coordinates and there will
be no ambiguities).

By a {\it time-dependent diffeomorphism} of a
manifold $U$ onto a manifold $Q$ we mean a diffeomorphism 
$
  \cC =(\cC_Q,\cC_\bR{}) : U\times\bR{} \to Q\times\bR{} 
$
whose second component $\cC_\bR{} : U\times\bR{}\to \bR{}$ is
the identity between the $\bR{}$-factors. The first component
$\cC_Q:U\times\bR{}$ is a differentiable map, that in the sequel we
denote $\cQ$. Thus,
$$
  \cC =(\cQ, \mathrm{id}_\bR{}) \;:\: U\times\bR{} \to Q\times\bR{}
  \,,
  \qquad
  \cC(u,t)  = (\cQ(u,t),t )
$$
and, for each $t$, the map
$$
   \cQ_t\;:=\; \cQ(\,\cdot\,,t)\;:\;U\to Q
$$
is a diffeomorphism (and ``smoothly depends on $t$'').

With the identifications $T(U\times\bR{}) \simeq TU\times T\bR{}$ and 
$T(Q\times\bR{}) \simeq TQ\times T\bR{}$, the tangent map
$T\cC: T(U\times\bR{}) \to T(Q\times\bR{})$ can be seen as a
diffeomorphism from $TU\times T\bR{}$ to $TQ\times T\bR{}$ whose
second component is the identity on the factor $T\bR{}$.
Restricting $T\cC$  to the
unit tangent vector in the $T\bR{}$-factor gives the diffeomorphism
$T\cC|_{\dot t =1}: TU\times \bR{}\times\{1\} \to
TQ\times \bR{}\times \{1\}$, that
we regard as a diffeomorphism
\begin{equation*}
  D\cC : TU\times \bR{} \to TQ\times \bR{} \,.
\end{equation*}
Explicitly, if for all $u\in U$ and $t\in\bR{}$ we write
$$
   \mathring\cQ_t(u):=\der{}t \cQ(u,t) \in T_{\cQ_t(u)}Q \,,
$$
then
\begin{equation*}
  D\cC(v_u,t) = 
  \big( \, T_u\cQ_t \cdot v_u + \mathring \cQ_t(u) ,\, t \, \big) 
  \in T_{\cQ_t(u)}Q\times\bR{} 
\end{equation*}
for all $u\in U$, $t\in \bR{}$ and $v_u\in T_uU$. Clearly,
$D\cC$ is a time-dependent diffeomorphism of $TU$ onto $TQ$.

In coordinates ($u\in U$, $q\in Q$) we will write $\cQ'$ for
$\der\cQ u$ and $\mathring\cQ$ for $\der\cQ t$ and, for given~$t$,
$\cQ'_t=\der\cQ u(\cdot ,t)$ and $\mathring\cQ_t=\der\cQ t(\cdot,
t)$. Thus, $\cQ'$ and $\mathring\cQ$ are defined on $U\times\bR{}$
while, for each $t$, $\cQ'_t$ and $\mathring\cQ_t$ are defined on
$U$. With this notation,
$
  T\cC(u,\dot u,t,\dot t) = 
  (\cQ(u,t),\, \cQ'_t(u)\dot u+\mathring\cQ_t(u)\dot t,\, t,\,
\dot t)
$
and
\begin{equation}\label{DC}
  D\cC(u,\dot u,t) = 
  (\cQ_t(u) ,\, \cQ'_t(u)\dot u+\mathring\cQ_t(u) ,\, t)  \,.
\end{equation}

If $\cC$ is a time-dependent diffeomorphism from $U$ onto $Q$, then the pull
back $\tilde M:=D\cC^{-1}(M)$ of an extended affine subbundle
$M=\{M_t\times\{t\}:t\in\bR{}\}$ of $TQ$ is an extended affine subbundle
of $TU$. In coordinates, if $M$ is described by $S(q,t)\dot q +
s(q,t)=0$ then $\tilde M$ is described by 
$
  \tilde S(u,t)\dot u + \tilde s(u,t)=0
$ 
with
\begin{equation}\label{Stilde}
  \tilde S = (S\circ\cC)\cQ' \,,\qquad
  \tilde s = s\circ \cC + (S\circ\cC) \mathring \cQ \,,
\end{equation}
as is verified requiring that, for each $t$, $(u,\dot u)
\in \tilde M_t$ if and only if $(\cQ_t(u),\cQ'_t(u)\dot
u+\mathring\cQ_t(u)) \in M_t$. 

The following fact is proven, in the time-independent case, in
\cite{FS2015}; the  generalization to the time-dependent case is
straightforward and we omit it.

\begin{proposition}\label{prop1} Consider a nonholonomic system
$X_{L,Q,M}$ and a time-dependent diffeomorphism $\cC$ from a manifold
$U$ onto $Q$. Then, the pull-back of $X_{L,Q,M}$ under the
restriction to $M$ of $D\cC$ coincides with the nonholonomic system 
$X_{\tilde L,U,\tilde M}$ with $\tilde L = L\circ D\cC$ and $\tilde
M=D\cC^{-1}(M)$.
\end{proposition}

\section{Conservation of moving energy }
\label{sec:ModifiedEnergy}

\subsection{Example}
We begin by illustrating the mechanism we have in mind on a sample
system---the well known sphere on a turntable considered by Pars
\cite{pars}, Neimark and Fufaev~\cite{NF} and several others, see
e.g. \cite{BKMM, cortes}.

This system is formed by a homogeneous solid sphere constrained to
roll without sliding on a table which, relatively to an inertial
reference frame, rotates with constant rate $\Omega$ around an axis
orthogonal to it. In the mentioned references, and in all the other
works we could find, the system is described with respect to the
inertial reference frame, but is nevertheless time-independent. Let
$\{O;x,y,z\}$ be such a frame, and assume that the table lies in the
$xy$-plane, and rotates about the $z$-axis. 

The configuration manifold $Q$ is $\bR2\times\mathrm{SO(3)}\ni(q,R)$
where $q=(x,y)$ are the coordinates of the point of contact between
sphere and table and $R$ is the attitude matrix of the sphere. We
identify $T\mathrm{SO(3)}$ and $\mathrm{SO(3)}\times\bR3$ via right
trivialization. The constraint manifold $M_0$ is 8-dimensional and is
diffeomorphic to $\bR2\times\mathrm{SO(3)}\times\bR3 \ni
(q,R,\omega)$, where $\omega=(\omega_x,\omega_y,\omega_z)$ is the
angular velocity in space of the sphere. Up to an inessential factor
$m$, the mass of the sphere, the Lagrangian is
\begin{equation}
\label{L-TT}
  L = \frac12\|\dot q\|^2 + \frac12 ca^2 \|\omega\|^2
\end{equation}
where $a$ is the radius of the sphere and $ca^2$, with
$c>0$, its moment of inertia (divided by $m$). The condition of rolling 
without sliding is given by the affine constraint
\begin{equation}
\label{eq:ball-table-constraint}
  \dot x = a \omega_y -\Omega y \,,\qquad 
  \dot y = -a \omega_x +   \Omega x \,.
\end{equation}

The Lagrangian and the constraints are
$\mathrm{SO(3)}$-invariant.\footnote{They are in fact invariant under 
an action of $\mathrm{SO(3)}\times S^1$, but this is not used in the
quoted references.} Reduction under this action consists merely in
cutting away the $\mathrm{SO(3)}$ factor and produces a
5-dimensional reduced system on $\bR{2}\times \bR{3} \ni
(q,\omega)$. The equations of motion of the reduced system are the
two equations (\ref{eq:ball-table-constraint}) and the three equations
$$
  \dot \omega_x = \frac \nu a \,(a\omega_y - \Omega y) \,,\qquad 
  \dot\omega_y = \frac \nu a\, (-a\omega_x + \Omega x) \,,\qquad
  \dot \omega_z= 0
$$
where $\nu = \frac{\Omega}{1+c}$, while the reconstruction equation is
$\dot\omega=\dot R R^T$. It is elementary to show that
the solutions of the reduced systems are periodic, with frequency
$\nu$. 

In some of the quoted references, e.g. in \cite{NF}, it is remarked
that the reduced equations have the three independent first integrals
$$
   \omega_z \,,\qquad 
   a\omega_x - \nu x   \,,\qquad
   a\omega_y - \nu y \,.
$$
However, the periodicity of a flow in a 5-dimensional phase space
implies the existence of four, not just three, independent first
integrals. To our knowledge (and surprise), this fact, and the
existence of a fourth independent first integral, do not seem to have
been noticed before.

The obvious candidate for the missing first integral would seem to be
the (projection to the reduced phase space of the) energy $E_{L,M_0}$
of the nonholonomic system, which is the restriction to the
constraint manifold $M_0$ of the energy $E_L$ of the Lagrangian $L$,
that in this case coincides with $L$. However, energy is not conserved: the
equations of motions give
$$
   \frac d{dt} E_{L,M_0} = ac\nu\Omega(x\omega_y - y\omega_x) \,.
$$
Nevertheless, a simple computation shows that the function 
\begin{equation}\label{fourth-integral}
  E_{L,M_0} - \Omega^2 (x^2+y^2) + \Omega a (x\omega_x+y\omega_y)
\end{equation}
is a first integral of the nonholonomic system. Being
$\mathrm{SO(3)}$-invariant, this function is also a first integral of
the reduced system, and it is independent of the previous three
(except where $\dot x=\dot y=0$). 

This additional first integral has a simple interpretation. In a
system of time-dependent, rotating coordinates in which the table is
at rest, the constraint of rolling without sliding on the table is
linear, and the Lagrangian, which is the pull back $\tilde L$ of
$L$, is still time-independent. Therefore, the energy is now
conserved and its push forward to the original coordinates is a first
integral, that turns out to coincide with (\ref{fourth-integral}). 

The reason why the push forward of the energy in the rotating
coordinates is different from the original energy (and may thus be
conserved) is due to the time-dependency of the coordinate change. The
Lagrangian in the rotating coordinates has the form $\tilde L=\tilde
L_2+\tilde L_1+\tilde L_0$, where the dependence of each $\tilde L_i$
on the velocities is homogeneous of degree~$i$.\footnote{Clearly,
$\tilde L_2$ may be interpreted as the kinetic energy, $-\tilde L_0$
as the potential energy of the centrifugal force and $-\tilde L_1$ as
the generalized potential of the Coriolis force in a rotating,
non-inertial reference frame in which the table is at rest. We prefer
changing coordinates, instead of reference frames, since this exempts
us from embedding the dependence on the choice of a reference frame
into the theory, as e.g. in \cite{pagani91}.} The function $\tilde
L_1$ does not contribute to the energy  $E_{\tilde L} = \tilde L_2-
\tilde L_0$ of $\tilde L$ and the push forward of $E_{\tilde L}$ to
the original coordinates differs from $E_L$ by the push forward of
$\tilde L_1$.

\begin{remark} {\rm A completely similar situation is met in the system
formed by a
vertical disk constrained to roll without sliding on a uniformly
rotating plane, considered in \cite{ferrario-passerini}. In that
reference, the system is actually studied in a rotating frame, where
the constraint is linear and the energy is conserved.  However, the
authors directly integrate the reduced equations of motion without
noticing the conservation of energy.
}\end{remark}

\subsection{Moving energy and its conservation} First recall that the
{energy} (or `Jacobi integral') of a Lagrangian $L:TQ\times\bR{}\to
\bR{}$ is the function $E_{L} : TQ\times \bR{} \to \bR{}$ given by
$$
  E_{L}(v_q,t) \,:=\, \langle p(v_q,t) , v_q\rangle_q - L(v_q,t) 
  \qquad \forall q\in Q \,,\; v_q\in T_qQ \,,\; t\in\bR{} \,,
$$
where $p$ is the momentum covector and $\langle\ ,\ \rangle_q$ denotes
the pairing between $T^*_qQ$ and $T_qQ$. In coordinates,  $E_{L} =
\dot q \cdot \der L {\dot q} - L$, where the dot denotes the scalar
product in $\bR n$. 

\begin{definition}
Let $X_{L,Q,M}$ be a (either time-dependent or
time-independent)
nonholonomic system with affine constraints.

\begin{list}{}
{\leftmargin2em\labelwidth1.2em\labelsep.5em\itemindent0em
\topsep0.5ex\itemsep-0.2ex}
\item[i.] The {\rm energy} $E_{L,M}:M\to \bR{}$ of $X_{L,Q,M}$ is the
restriction of $E_L$ to $M$:
\begin{equation*}
  E_{L,M} := E_L|_M  \,.
\end{equation*}

\item[ii.] If $\cC:U\times \bR{} \to Q\times \bR{}$ is a
time-dependent diffeomorphism, then the {\rm moving energy} of $X_{L,Q,M}$
induced by $\cC$ is the restriction $E^*_{L,\cC,M}$ to $M$ of the function
\begin{equation}\label{moven}
   E^*_{L,\cC} := E_{L\circ D\cC}\circ D\cC^{-1} 
   : TQ\times\bR{} \to \bR{} \,.
\end{equation}
\end{list}
\end{definition}

\begin{proposition}\label{prop2} In the hypotheses of the above
definition,
\begin{equation}\label{DiffIntJac}
  E^*_{L,\cC} = E_L
  - 
  \langle p, \mathring\cQ \circ\cC^{-1} \rangle \,.
\end{equation}
\end{proposition}

\begin{proof} The proof can be done in coordinates. If
$\tilde L=L\circ D\cC$, then, from (\ref{DC}), 
$\tilde L(u,\dot u,t) = L(\cQ_t(u),\cQ_t'(u)\dot u + \mathring\cQ_t(u),t)$. 
Thus
$
  E_{\tilde L} (u,\dot u,t) 
  = 
  \dot u\cdot \der {\tilde L}{\dot u} (u,\dot u,t)
  - 
  \tilde L(u,\dot u,t)
  =
  \cQ'_t(u) \dot u \cdot \der L {\dot q}(D\cC(u,\dot u,t)) 
  - 
  L(D\cC(u,\dot u,t)) 
  =
  E_L(D\cC(u,\dot u,t)) 
  - 
  \mathring\cQ_t(u)\cdot \der L {\dot q} (D\cC(u,\dot u,t)) 
$.
In (\ref{DiffIntJac}) we have written $\mathring\cQ \circ\cC^{-1}$
instead of $\mathring\cQ \circ D\cC^{-1}$ because $\cQ$ is
independent of the velocities. 
\end{proof}

The interest of considering a moving energy $E^*_{L,\cC,M}$
resides in the fact that the function $E^*_{L,\cC}$
differs from the energy $E_L$ of the
Lagrangian $L$ by a term which is produced by the time-dependence of
the diffeomorphism.\footnote{A fact which is well known in the theory
of time-dependent canonical transformations.} It is therefore
possible that the function $E^*_{L,\cC,M}$ is a first integral even
if $E_{L,M}$ is not. We now formalize this possibility in the case
of {time-independent} nonholonomic systems:

\begin{theorem}
\label{thm1}
Consider a time-independent nonholonomic system with affine constraints
$X_{L,Q,M_0\times\bR{}}$ and a time-dependent diffeomorphism $\cC$ from a
manifold $U$ to $Q$. Assume that:
\begin{list}{}
{\leftmargin2em\labelwidth1.2em\labelsep.5em\itemindent0em
\topsep0.5ex\itemsep-0.2ex}
\item[i.] $L\circ D\cC$ is independent of $t$.
\item[ii.] $E^*_{L,\cC}$ is independent of $t$.
\item[iii.] $D\cC^{-1}(M_0\times \bR{})$ is a time-independent
extended linear subbundle of $TU$. 
\end{list}
Then, the moving energy $E^*_{L,\cC,M_0\times\bR{}}$  
is a time-independent first integral of $X_{L,Q,M_0\times\bR{}}$.
\end{theorem}

\begin{proof} 
Hypothesis iii{.} means that
$D\cC^{-1}(M_0\times \bR{})=\tilde M_0\times \bR{}$
for a fixed linear subbundle $\tilde M_0$ of $TU$.  By Proposition
\ref{prop1}, $D\cC$ conjugates $X_{L,Q,M_0\times\bR{}}$ to the
nonholonomic system  $X_{\tilde L, U, \tilde M_0\times\bR{}}$ with
$\tilde L = L\circ D\cC$. By hypotheses i. and iii.,
$X_{\tilde L, U, \tilde M_0\times\bR{}}$
has linear constraints and is time-independent.
Therefore, the energy $E_{\tilde L,\tilde M_0\times\bR{} }$ is  a
(time-independent) first integral of $X_{\tilde L, U,\tilde M_0\times\bR{}}$. It
follows that its push-forward
$E_{\tilde L,\tilde M_0\times\bR{} }\circ (D\cC^{-1}|_{M_0\times\bR{}})$ 
is a first integral of
$X_{L,Q,M_0\times\bR{}}$. Since $D\cC$ maps $\tilde M_0\times\bR{}$
diffeomorphically onto $M_0\times\bR{}$,
$$
  E_{\tilde L,\tilde M_0\times\bR{} }\circ (D\cC^{-1}|_{M_0\times\bR{}})
  =
  (E_{\tilde L,\tilde M_0\times\bR{} }\circ D\cC^{-1})|_{M_0\times\bR{}}
  =
  E^*_{L,\cC,M_0\times\bR{}} \,.
$$ 
This proves that $E^*_{L,\cC,M_0\times\bR{}}$ is a first integral
of $X_{L,Q,M_0\times\bR{}}$. Hypothesis ii. ensures that it is a
time-independent function.  \end{proof}

We have stated Theorem 1 in terms of time-independent nonholonomic
systems on the extended phase space so as to properly regard, in
hypotheses i. and ii., the functions $L$ and $E^*_{L,\cC}$ as defined
on $TQ\times\bR{}$, even though constant on $\bR{}$. But if we
identify functions on $TQ$ and functions on $TQ\times\bR{}$ that are
constant on $\bR{}$, then Theorem 1 states that, {\it under
hypotheses i., ii. and iii., a time-independent
nonholonomic system with affine
constrains $X_{L,Q,M_0}$ has the time-independent first integral 
$E^*_{L,\cC,M_0}$} (in fact, since $E^*_{L,\cC}$ is time-independent,
$E^*_{L,\cC,M_0\times\bR{}}$ and $E^*_{L,\cC,M_0}$  may be
identified). From now on, we will adopt this point of view.

\begin{remarks}
{\rm  
(i) In time-dependent nonholonomic systems, either with linear or with
affine constraints, the energy $E_{L,M}$ is ordinarily
time-dependent; even though it is not impossible that it is a
(time-dependent) first integral, we do not consider this case because,
in our opinion, from  a dynamical point of view only time-independent
first integrals are of interest.

(ii) One might weaken the hypotheses of Theorem \ref{thm1} in various
ways, e.g. by requiring that only the restriction of $E^*_{L,\cC}$ to
$M_0\times \bR{}$  be time-independent. However, the setting of Theorem
\ref{thm1} is sufficient for our application in Section
\ref{sec:ball-turning-surface}
}\end{remarks}

\subsection{On the conditions of Theorem \ref{thm1}} 

The situation of Theorem \ref{thm1} might appear very special. Our
next goal, in Section 3.4, is to show that such a
situation is instead easily verified in presence of
symmetries. In order to gain some insight on this possibility, we begin
by establishing conditions under which $E^*_{L,\cC}$ is
time-independent and conditions under which the extended constraint
submanifold $D\cC^{-1}(M_0\times\bR{})$ in the new coordinates is
linear (even though possibly time-dependent).

\begin{proposition}\label{prop3} 
Consider a time-independent nonholonomic system $X_{L,Q,M_0}$ with
affine constraints. Denote by $\xi_0+\cD_0$ the affine distribution
on $Q$ that corresponds to $M_0$, with $\xi_0$ a vector field on $Q$
and $\cD_0$ a distribution on $Q$. Consider a time-dependent
diffeomorphism $\cC=(\cQ,\mathrm{id}_\bR{})$ of $Q$ to itself.

\begin{list}{}
{\leftmargin2em\labelwidth1.2em\labelsep.5em\itemindent0em
\topsep0.5ex\itemsep-0.2ex}

\item[i.] $E^*_{L,\cC}$ is time-independent if and only if $\cQ$
is the flow of a vector field $Y$ on $Q$, and in that case
$$
  E^*_{L,\cC} = E_L - \langle p,Y \rangle \,.
$$

\item[ii.] Assume that $\cQ$ is the flow of a vector field
$Y$ on $Q$. Then, $D\cC^{-1}(M_0\times\bR{})$ is an extended linear
subbundle of $TQ$ if and only if the vector field $Y-\xi_0$ is a
section of~$\cD_0$.

\end{list}
\end{proposition}

\begin{proof} 
(i) First note that
$
  \mathring\cQ\circ \cC^{-1}(q,t) = 
  \mathring\cQ_t\circ\cQ_t^{-1}(q) 
$
for all $q$ and $t$.  Let us write $Y(q,t) =
\mathring\cQ_t\circ\cQ_t^{-1}(q)$ and note that $Y(q,t)\in T_qQ$.
Since $L$ and $E_L$ are time-independent, it follows from 
(\ref{DiffIntJac}) that $E^*_{L,\cC}$ is time-independent if and
only if $\langle Y,p\rangle$ is time-independent, that is, given that
the momentum covector $p$ does not depend on time, if and only if 
$$
  \der{}t \big\langle Y , p \big\rangle 
  =
  \Big\langle \der{Y}t , p \Big\rangle \,:\, 
  TQ\times\bR{} \to \bR{} 
$$
vanishes. Since the Lagrangian $L$ is regular, the map
$v_q \mapsto \langle p,v_q\rangle_q$ is a local
diffeomorphism for each $q$. Therefore, $\langle \der{Y}t , p
\rangle$  vanishes identically in $TQ\times\bR{}$  if and
only if $\der Yt=0$. This shows that the time-independence of
$E^*_{L,\cC}$ is equivalent to $\mathring\cQ_t\circ\cQ_t^{-1}=Y$ with
$Y$ independent of $t$. But then $Y$ is a vector field on $Q$
and, since $\mathring\cQ=Y\circ\cQ$, $\cQ$
is the flow of $Y$.

(ii) The proof can be done in coordinates. 
Let $M_0$ be given by $S(q)\dot q+s(q)=0$. Then 
$D\cC^{-1}(M_0\times\bR{})$ is an affine subbundle of $TQ$ that is
described by $\tilde S(u,t)\dot
u+\tilde s(u,t)=0$ with $\tilde S$ and $\tilde s$ as in 
(\ref{Stilde}). Its linearity is equivalent to the vanishing of
$\tilde s= s\circ\cC + (S\circ\cC)\mathring\cQ$, that is, given that
$S$ and $s$ are time-independent, $s\circ\cC =s\circ\cQ $ and
$S\circ\cC=S\circ\cQ$, to the vanishing of
$s + S (\mathring \cQ_t\circ\cQ_t^{-1}) = S(Y-\xi_0)$.
\end{proof}

This proposition suggests that, in order to obtain a time-independent
conserved moving energy, the time-dependent diffeomorphism $\cC$
should be constructed as the flow of a vector field on $Q$ that
differs from the vector field $\xi_0$ by a section of the
distribution $\cD_0$. The freedom in the choice of this section might
then be used to try to make $L\circ D\cC$ and
$D\cC^{-1}(M_0\times\bR{})$ time-independent.  In the next section we
will show that this is always possible if the system admits a
symmetry group, with suitable properties, by choosing  $Y$ as an
infinitesimal generator of the group action, that is, by choosing
$\cQ$ as the action of a one-parameter subgroup.

\subsection{Symmetry and conservation of moving energy}

We consider now a time-independent nonholonomic system $X_{L,Q,M_0}$
with affine constraints whose Lagrangian and constraint distribution
are invariant---in a sense made precise in Hypotheses (H1) and (H2)
below---under an action $\Psi:G\times Q \to Q$ of a Lie group $G$ on
$Q$. As in Proposition \ref{prop3}, we denote by $\cM_0=\xi_0 +
\cD_0$ the affine distribution on $Q$ corresponding to~$M_0$.

For each $q\in Q$, we write as usual $\Psi_g(q)$ for $\Psi(g,q)$. We
denote by $\Psi^{TQ}:G\times TQ\to TQ$ the tangent lift of the action
$\Psi$, which is the action of $G$ on $TQ$ given by
$$
   \Psi^{TQ}_g(v_q) = T_q\Psi_g \cdot v_q 
$$
(in coordinates, $\Psi^{TQ}_g(q,\dot q) = \big( \Psi_g(q),
\Psi'_g(q)\dot q\big)$ with $\Psi'_g = \der{\Psi_g}q$).  We make the
following two hypotheses:
\begin{list}{}
{\leftmargin2.5em\labelwidth2.1em\labelsep.5em\itemindent0em
\topsep0.5ex\itemsep-0.2ex}

\item[(H1)] $L$ is invariant under $\Psi^{TQ}$, namely
$$
    L\circ \Psi^{TQ}_g = L \qquad \forall g\in G
$$
(in coordinates, $L ( \Psi_g(q) , \Psi'_g(q)\dot q ) = L(q,\dot q)$ $\forall
g, q,\dot q$). 

\item[(H2)] The distribution $\cD_0$ is invariant under $\Psi$, in
the sense that
\begin{equation*}
  (\cD_0)_{\Psi_g(q)} = T_q\Psi_g \cdot (\cD_0)_q \qquad \forall g\in G \,, q\in Q
\end{equation*}
(we need not make any hypothesis on the nonhomogeneous term $\xi_0$ and
on the invariance of $M_0$ under the group action).

\end{list}
Under these hypotheses, it is rather natural to try to build the
time-dependent
diffeomorphism $\cC=(\cQ,\mathrm{id}_\bR{})$ that leads to a
conserved moving energy by choosing $\cQ$ as a one-parameter subgroup of
the action~$\Psi$. 

For $\eta\in\mathfrak g$, the Lie algebra of
$G$, denote by
$$
    Y_\eta := \frac d{dt}\Psi_{\exp(t\eta)}|_{t=0}
$$ 
the infinitesimal generator of the action of the one-parameter
subgroup generated by $\eta$ and by
$$
  J_\eta := \langle p , Y_\eta \rangle 
$$
the momentum map of the lifted action of the same one-parameter
subgroup. The moving energy of $X_{L,Q,M_0}$ relative of the
time-dependent diffeomorphism
$\cC_\eta=(\Phi^{Y_\eta},\mathrm{id}_\bR{})$, where 
$\Phi^{Y_\eta}:Q\times\bR{}\to Q$ is the flow of $Y_\eta$, is thus
the restriction to $M_0$ of the function
\begin{equation}\label{En-e-MomMap}
  E^*_{L,\cC_\eta} = E_L - J_\eta \,.
\end{equation}

\begin{theorem}\label{thm2}
Consider a time-independent nonholonomic system $X_{L,Q,M_0}$ with affine
constraints and an action $\Psi$ of a Lie group
$G$ on $Q$. Assume (H1), (H2) and 
\begin{list}{}
{\leftmargin2.5em\labelwidth2.1em\labelsep.5em\itemindent0em
\topsep1.5ex\itemsep-0.1ex}

\item[(H3)] $\eta\in\mathfrak g$ is such that $Y_\eta-\xi_0$ is a
section of $\cD_0$.

\end{list}
Then, the moving energy $ E^*_{L,\cC_\eta,M_0}$
is a time-independent first integral of $X_{L,Q,M_0}$.
\end{theorem}

\begin{proof} Let $\cQ$ be the flow of $Y_\eta$. The conclusion
follows if we show that the three hypotheses of Theorem \ref{thm1}
are satisfied with $\cC=(\cQ,\mathrm{id}_\bR{})$. We may check them
in coordinates.

{\it Hypothesis i.} With this choice of $\cQ$, $\cQ_t =
\Psi_{\exp(t\eta)}$ for all
$t$. Hence, by hypothesis (H1), $L(\cQ_t(q),\cQ'_t(q)\dot q) = L(q,\dot
q)$ for all $q,\dot q,t$. Define $\tilde L=L\circ D \cC$. Since a
vector field is invariant under its own flow, $\mathring\cQ_t =
Y_\eta\circ \cQ_t  = \cQ'_t Y_\eta$. Thus
$$
  \tilde L( q,\dot q , t) 
  = 
  L\big( \cQ_t(q) , \cQ'_t(q)\dot q + \mathring\cQ_t(q) \big)
  = 
  L\big( \cQ_t(q) , \cQ'_t(q)[ \dot q + Y_\eta(q)] \big)
  = 
  L( q , \dot q + Y_\eta(q))
$$
that shows that $\tilde L$ is time-independent. 

{\it Hypothesis ii.} Under hypothesis (H3), by item {ii.} of Proposition
\ref{prop3},  $\tilde M = D\cC^{-1}(M_0\times\bR{})$ is an extended
linear subbundle of $TQ$. We prove that this subbundle is
time-independent.
If $M_0$ is described by $S(q)\dot q+ s(q)=0$,
then its associated distribution $\cD_0$ has fibers $\ker S(q)$
and hypothesis (H2) is
$\ker S(\Psi_g(q)) = \Psi'_g(q)\,\ker S(q)$ for all $g$ and $q$, or
\begin{equation}\label{InvD2}
  \ker \big[ S(\cQ_t(q)) \big]
  = 
 \cQ'_t(q)\,\ker [S(q)] \qquad \forall t\,, q \,.
\end{equation}
In turn, the extended distribution $\tilde\cM$ associated to $\tilde M$
has fibers $(\tilde \cM_t)_q=\ker[ \tilde S(q,t)]$ with 
$$
  \tilde S (q,t) = S(\cQ_t(q))\cQ'_t(q) \,,
$$
see (\ref{Stilde}). Hence, using (\ref{InvD2}) and the fact that, if
$B$ is an invertible $n\times n$ matrix and $S$ is a $k\times n$
matrix, then $\ker(SB) =B^{-1} \ker S$, we have
$$
  \ker\tilde S(q,t) 
  \;=\; 
  \ker\big[ S(\cQ_t(q))\cQ'_t(q)\big] 
  \;=\; 
  \cQ'_t(q)^{-1} \ker S(\cQ_t(q)) 
  \;=\; 
  \ker S(q) \,.
$$
Thus $\tilde\cM_t=\cD_0$ for all $t$ and $\tilde M$ is
time-independent. 

{\it Hypothesis iii.} This follows from item i. of Proposition
\ref{prop3}.
\end{proof}

\begin{remark}{\rm 
The condition that $Y_\eta-\xi_0$ is a section of $\cD_0$ implies
that the orbits of the group action $\Psi^{TQ}$ must be transversal
to the constraint manifold.
}
\end{remark}

\begin{example} {\rm
The sphere on the turntable of Section~3.1 is an instance of the 
situation of Theorem 2. As in that section, we identify the tangent
spaces to $Q$ with $\bR5$ via right-trivialization of
$T\mathrm{SO(3)}$. The natural symmetry group of the problem is
$S^1\times \mathrm{SO(3)}$ (that acts as in (\ref{azione}) below), but
for the sake of applying Theorem 2 we may consider only its subgroup
$G=S^1\times\{\mathbb I\}$, that acts on $Q = \mathbb{R}^{2} \times
\textrm{SO}(3)$ as
$$
  \Psi_{\theta} (q, R) = (H_{\theta} q, H_\theta R) \,.
$$
Here $H_\theta$ is the matrix of the anticlockwise rotation of angle
$\theta$ around the $z$-axis and, with a small abuse, we
identify vectors $(x,y)\in\bR2$ with vectors $(x,y,0)\in\bR3$. The
(right-trivialized) infinitesimal generator of the action that
corresponds to the Lie algebra element $\eta \in\bR{}$ is $Y_\eta =
(-\eta y,\eta x,0,0,\eta )$, and the corresponding momentum is $J_\eta
=\eta (x\dot
y-y\dot x+ca^2\omega_z)$.  The (right-trivialized) tangent lift of
this action is
$$
  \Psi^{TQ}_{\theta}(q,R,\dot q,\omega) 
  =
  (H_{\theta} q, H_\theta R, H_{\theta} \dot q, H_\theta \omega ) 
$$
and leaves the Lagrangian (\ref{L-TT}) invariant, as in (H1). The
distribution $\cD_0$ associated to the affine constraint
(\ref{eq:ball-table-constraint}) is given by $\dot x = a \omega_y$,
$\dot y=-a\omega_x$ (that is, $\dot q =a\omega\times e_z$
with $\times$ the vector product in $\bR3$) and is invariant under the
action $\Psi_\theta$, as in (H2). Finally, the
nonhomogeneous part of the constraint (\ref{eq:ball-table-constraint})
is the vector field $\xi_0 = (-\Omega y, \Omega x, 0, 0, 0)$
and $Y_\Omega - \xi_0=(0,0,0,0,\Omega)$ lies
in $\cD_0$, as in (H3). By Theorem~2, the moving energy
$E^*_{L,\cC_\Omega,M_0} = E_{L,M_0} - J_\Omega|_{M_0}$ is
conserved. Using (\ref{eq:ball-table-constraint}) one verifies that
this moving energy coincides with~(\ref{fourth-integral}),
up to a constant term $ca^2\Omega\omega_z$.
}\end{example}

\subsection{Connection to (the nonholonomic failure of) Noether theorem}

In the setting of Theorem~\ref{thm2} it is natural to view the
conservation of the moving energy
$$
   E^*_{L,\cC_\eta,M_0}= E_{L,M_0} - J_{\eta,M_0} \,,
$$
where 
$$
  J_{\eta, M_0} = J_\eta|_{M_0} \,,
$$
as related to the invariance of the system under the action $\Psi$. 
If the energy $E_{L,M_0}$ is not conserved, then the conservation of
the moving energy $E^*_{L,\cC_\eta,M_0}$ is only possible if 
$J_{\eta,M_0}$ is not conserved. Thus, Theorem \ref{thm2} produces
a conserved quantity from the sum of two quantities---the energy and
a component of the momentum map---that, at variance
from what would happen if the system was holonomic, are not conserved.
In a way, this first integral seems to be produced
notwithstanding---or perhaps thanks to---the failure in nonholonomic
mechanics of two cornerstones of Lagrangian mechanics:  conservation
of energy and Noether theorem.

It has some interest to understand why, at least in the symmetric case
considered here, the mechanism of Theorem  \ref{thm1}, that obviously has
no interest in the Lagrangian case, is instead of interest in the
nonholonomic case. In the absence of nonholonomic constraints
($M_0=\tilde M_0=TQ$), Theorem \ref{thm1} is obviously true: the
time-independence of the two Lagrangians $L$ and $\tilde L$ implies
that both functions $E_L$ and $E^*_{L,\cC_\eta}$ are first integrals of the
Lagrangian system described by the Lagrangian $L$. However, since 
in Lagrangian mechanics the momentum map of a lifted action that
leaves the Lagrangian invariant is conserved, this mechanism
can hardly be seen as disclosing a `new' first
integral $E^*_{L,\cC_\eta}$: the difference $E^*_{L,\cC_\eta}-E_L$
is a component of the momentum
map, and it is thus a first integral for the very same reason of
symmetry that underpins the possibility of passing to moving
coordinates without introducing a time-dependence in the Lagrangian.

Explaining why things are different in the nonholonomic case---and how
they are different---requires exploiting the role of the reaction
forces, along the lines of \cite{FRS,FS2015}. The hypothesis of
ideality of the constraints assumes that the constraint can---a
priori---exert all reaction forces that lie in $(\cD_0)_q^\circ$.
However, the set of reaction forces that is actually exerted when the
system $X_{L,Q,M_0}$ is in a configuration $q\in Q$ with any
possible velocity $\dot q\in (\cM_0)_q$ is given by
\[
  \cR_{q}:= 
  \bigcup_{\dot q\in (\cM_0)_q} R_{L,M}(q,\dot q) 
\]
and can be (and typically is) smaller than $(\cD_0)_q^\circ$. Therefore, the
annihilators of $\cR_{q}$ can be (and typically are) larger than
the fibers of $(\cD_t)_q$. These annihilators are the fibers of a
distribution $\cR^\circ$ on $Q$, which was introduced in \cite{FRS} (in
the case of time-independent linear constraints, but the generalization
to the case of time-independent affine constraints is straightforward
\cite{FS2015}) and was called the {\rm
reaction-annihilator distribution}. We refer to these works for
further details and we limit ourselves to note that
$$
  (\cD_0)_q  \subseteq \cR_{q}^\circ \qquad \forall \, q,t \,.
$$

\begin{proposition}\label{prop4} {\rm \cite{FS2015}} Consider a
time-independent nonholonomic system with affine constraints
$X_{L,Q,M_0}$, and let $\cM_0=\cD_0+\xi_0$. 
\begin{list}{}
{\leftmargin2em\labelwidth1.2em\labelsep.5em\itemindent0em
\topsep0.5ex\itemsep-0.2ex}
\item[i.] The energy $E_{L,M_0}$ is conserved if and only if $\xi_0$ is a section of
$\cR^\circ$.
\item[ii.] Assume that the Lagrangian $L:TQ\to\bR{}$ is invariant under the
tangent lift of an action of a Lie group $G$ on $Q$, namely
$L\circ\Psi^{TG}_g=L$ for all $g\in G$. Then, for any
$\eta\in\mathfrak g$, $J_{\eta,M_0}$ is a first
integral of $X_{L,Q,M_0}$ if and only if $Y_\eta$ is a section of
$\cR^\circ$.
\end{list}
\end{proposition}

Assume, thus, that the energy of a time-independent nonholonomic system
$X_{L,Q,M_0}$ with affine constraints is not conserved. By
Proposition \ref{prop4}, this happens if and only if $\xi_0$ is not a
section of $\cR^\circ$. By Theorem \ref{thm2}, under Hypotheses
(H1) and (H2), the existence of a conserved moving energy
$E^*_{L,\cC_\eta,M_0}$ depends
on  the existence of an infinitesimal generator $Y_\eta$ such that
the difference $\xi_0-Y_\eta$ is a section of $\cD_0$. Since
$\cD_0\subseteq\cR^\circ$ and $\xi_0$ is not a section of $\cR^\circ$,
this necessarily requires that $Y_\eta$ has a nonzero component off
$\cR^\circ$ and, still by Proposition \ref{prop4}, the momentum
$J_{\eta,M_0}$ is not conserved. Thus, at the basis of
the non-conservation of the energy and of the component of the
momentum map---that makes it possible for them to add up to give a
conserved function---there is the same reason: $\xi_0$ is not a
section of $\cR^\circ$.

\begin{remark} {\rm One might take as well a different point of view, and
see the function $E^*_{L,\cC_\eta}$
as the momentum map of the action, in the extended phase space,
given by a combination of
time-translation and of the lift of a one-parameter subgroup of $\Psi$. 
}\end{remark}

\section{Integrability of a sphere rolling on a rotating surface of revolution}
\label{sec:ball-turning-surface}
We outline now an application of the existence of a conserved moving
energy to the class of
systems, considered by Borisov, Mamaev and Kilin in \cite{BMK2002},
that are formed by a heavy homogeneous solid sphere constrained to rotate
without sliding on a moving surface of revolution; specifically, the
surface rotates---relatively to an inertial frame---with uniform
angular velocity $\Omega$ around its figure axis, which is assumed to
be vertical (that is, directed like gravity). Describing the system
in an inertial frame $\{O;x,y,z\}$ and using
time-independent coordinates, as done in \cite{BMK2002}, leads to a
time-independent Lagrangian; we will assume that the $z$-axis
coincides with the figure axis of the surface, see Figure 1.
The case $\Omega=0$ is classical. Its study goes back to
Routh \cite{routh}, who also considered special cases with
$\Omega\not=0$ (see also \cite{NF}). When
$\Omega\not=0$ and the surface is a horizontal plane, the system
reduces to the sphere on the turntable of section 3.1.

In all cases, either with $\Omega=0$ or with $\Omega\not=0$,
the (time-independent) constraint manifold $M_0$ has dimension~8 and the
system has an obvious $S^1\times\mathrm{SO(3)}$ symmetry. Reduction
leads to a four-dimensional system which has two independent first
integrals, that we will denote $J_1$ and $J_2$.
We will use the same symbols $J_1$ and $J_2$ also for the lifts of
these functions to the unreduced phase space, which are first
integrals of the unreduced system.
The existence of the two first integrals $J_1$ and $J_2$, when
$\Omega=0$, was known already to Routh \cite{routh,NF}. Their
existence when $\Omega\not=0$ has been proven in~\cite{BMK2002}. 

\begin{figure}
\vskip 5truecm
\begin{picture}(0,0)
{\small
\put(45,0){\scalebox{.85}{\includegraphics*{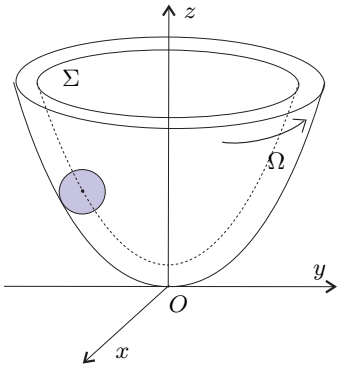}}}
\put(53,37){$\Sigma$}
\put(80,26){$\Omega$}
\put(60,1){$x$}
\put(86,12){$y$}
\put(69,46){$z$}
\put(67,7){$O$}
}
\end{picture}
\caption{\small The sphere in a rotating convex surface of revolution.}
\end{figure}

The integrability of the case $\Omega=0$ has been studied by Hermans
\cite{hermans} and Zenkov \cite{zenkov}. When $\Omega=0$ the
constraint is linear and the energy $E_{L,M_0}$ is conserved. As a
result, the unreduced system has the three first integrals $J_1$,
$J_2$ and $E_{L,M_0}$, which are independent in an open subset $M_0^*$ of
the phase space $M_0$ (specifically, they are independent everywhere
except on motions in which the center of the sphere either moves
horizontally or is at rest). Since the function $E_{L,M_0}$ is
$S^1\times\mathrm{SO(3)}$-invariant,  the reduced system has three
independent integrals as well.  Moreover, if the surface is upward
convex and the sphere rolls inside it, as in Figure~1, then the common
level sets of these three integrals in the four-dimensional reduced
phase space are compact, and hence are closed curves, and the reduced
dynamics is periodic \cite{hermans,zenkov} (reference \cite{hermans}
uses a different argument; for details on the properties of
the first integrals see 
\cite{FGS2005}). Since the symmetry group is compact and acts freely
on $M_0^*$, this in turn implies that the unreduced dynamics in $M_0^*$
is quasi-periodic on tori of dimension up to three. This was proven
in \cite{hermans} using a  reconstruction result from periodic
dynamics, originally due to Field and Krupa (see particularly
\cite{krupa,field,hermans,FG2007,CDS,field-book}). Reference \cite{zenkov} reaches the
same conclusion, but restricted to the motion of the center of mass,
that undergoes quasi-periodic motions on tori of dimension up to two.

When $\Omega\not=0$ the constraint is affine, not linear, and even if
the constraint and the Lagrangian are time-independent the
energy is not conserved. Therefore, the argument used for the case
$\Omega=0$ is not directly applicable. A different approach has been
taken by Borisov, Mamaev and Kilin, who proved that the reduced system
has an invariant measure and, using Jacobi theorem \cite{AKN-3d-ed},
deduced from this and from the existence
of the two first integrals $J_1$ and $J_2$ that the reduced system is 
integrable by quadratures \cite{BMK2002}. 

But as we now prove, a conserved `moving energy'
exists in this problem, and much stronger results can be obtained. 
Leaving for a future work a detailed study of the problem,
we limit here ourselves to some conclusions that
can be obtained combining Theorem \ref{thm1} 
with some general
arguments (essentially, continuity from the case $\Omega=0$):

\begin{theorem}\label{thm3} Consider the system formed by a heavy
homogeneous solid sphere that rolls without sliding on a surface of revolution,
which rotates with constant angular velocity $\Omega$ around its
figure axis, aligned with gravity. Then,
at least for $\Omega$ not too large:
\begin{list}{}
{\leftmargin2em\labelwidth1.2em\labelsep.5em\itemindent0em
\topsep0.5ex\itemsep-0.2ex}
\item[1.] The reduced system has three first integrals, which are
independent in some open nonempty subset of the four-dimensional
reduced phase space.
\end{list}
If, moreover, the surface is upward convex, and the sphere rolls
inside it, then: 
\begin{list}{}
{\leftmargin2em\labelwidth1.2em\labelsep.5em\itemindent0em
\topsep0.5ex\itemsep-0.2ex}
\item[2.] There is a nonempty open subset of the reduced phase space
where the reduced dynamics is periodic.
\item[3.] There is a nonempty open subset of the phase space of the
unreduced system in which motions are quasi-periodic, on tori of dimension
up to three. 
\end{list}
\end{theorem} 

\begin{proof}
Let $r\in\bR 3$ be the vector of the coordinates of the center of the
sphere relative to the considered inertial frame $\{O;x,y,z\}$, 
$R\in \mathrm{SO(3)}$ be the matrix that fixes the attitude of the
sphere relatively to that frame, and $\omega\in\bR 3$ be the 
angular velocity in space of the sphere relative to that frame.

The holonomic constraint that the sphere is in contact with the
surface of revolution can be modelled by imposing that the vector $r$
is constrained to a (fixed, time independent) surface of revolution
$\Sigma$, that we embed in $\bR3\ni r$. The configuration
space of the holonomic system is thus $Q=\Sigma\times\mathrm{SO(3)}
\ni (r,R)$ and the phase space $TQ$ can be identified with
$\Sigma\times\bR2\times \mathrm{SO(3)}\times\bR 3 \ni (r,\dot
r,R,\omega)$. The Lagrangian is the restriction to
$\Sigma\times\bR2\times \mathrm{SO(3)}\times\bR 3$ of the function
$L:\bR3\times\bR3\times \mathrm{SO(3)}\times\bR 3\to\bR{}$ given by
$$
  L(r,\dot r,R,\omega) 
  = 
  \frac 12 \|\dot r\|^2 + \frac 12 ca^2\|\omega\|^2 - g r_3 \,,
$$  
where the constants have obvious meanings, see (\ref{L-TT}), and is
time-independent. 

The constraint of rolling without sliding leads to a time-independent
nonholonomic system with affine constraints, with constraint
submanifold an 8-dimensional
affine subbundle $M_0$ of $\Sigma\times\bR2\times
\mathrm{SO(3)}\times\bR 3$ and Lagrangian the restriction of $L$ to
$M_0$. The Lagrangian and the constraint manifold are invariant under
the tangent lift of the action of $S^1\times\mathrm{SO(3)}$ on $Q$
given by 
\begin{equation}\label{azione}
  \Psi_{(\theta,P)}(r,R) = (H_\theta r ,H_\theta RP) 
\end{equation}
where $H_\theta$ is the $3\times 3$ matrix of rotation of angle
$\theta$ around the third axis. Once the kinematical states with the
sphere sitting at the point $r=0$ and spinning about the $z$-axis have
been removed from phase space, to prevent the need for singular
reduction, the (regularly) reduced phase space is a 4-dimensional
submanifold of $\Sigma\times\bR2\times\bR3\ni(r,\dot r,\omega)$.

We pass now to time-dependent coordinates $(s,S)$ in $Q$ with
\begin{equation}\label{cc}
  \cQ(s,S,t) = \big(H_{\Omega t}s ,  H_{\Omega t}S \big) 
\end{equation}
and lift them to a time-dependent coordinate change 
$D\cC:(s,\dot s,S,\nu,t) \mapsto (r,\dot r,R,\omega,t)$ in $TQ$.
In this coordinate system the surface is at rest; therefore, the
constraint of rolling without sliding is linear and time-independent
and defines a linear subbundle $\tilde M_0$ of $\Sigma\times\bR2\times
\mathrm{SO(3)}\times\bR 3$. Due to the symmetry of the system, the
Lagrangian $\tilde L =L\circ D\cC$ is time-independent as well.
(Incidentally, since $\cQ$ is the restriction of the action $\Psi$ to
a one-parameter subgroup, the time-independence of $\tilde L$ follows
as well from the argument used in the proof of Theorem \ref{thm2};
we mention also that 
$$
  \tilde L(s,\dot s,S,\nu) 
  = 
  \frac 12 \|\dot s+\Omega e_z\times s\|^2 
  + 
  \frac 12 ca^2\|\nu+\Omega e_z\|^2 - g s_3 \;\big).
$$  
Thus,
$X_{\tilde L, Q ,\tilde M_0}$ is a time-independent nonholonomic
system with linear constraints. Moreover, by item i. of Proposition
\ref{prop3} the function  $E_{L,\cC}^*$ is time-independent. By
Theorem~\ref{thm1}, we conclude that the original system
$X_{L,Q,M_0}$ has the conserved moving energy $E^*_{L,\cC,M_0}$. 

We may now prove the three statements.

1. Since the constraint manifolds $M_0$ and $\tilde M_0$, the Lagrangians
$L$ and  $\tilde L$ and the change of coordinates are $S^1\times
\mathrm{SO(3)}$-invariant, the function $E^*_{L,\cC,M_0}$ has this
very same
invariance property and descends to a first integral of the reduced
system. The reduced system has therefore the three first integrals
$J_1$, $J_2$ and $E^*_{L,\cC,M_0}$.
From the expressions for $J_1$ and $J_2$
given in \cite{BMK2002} and from the expression above of $\tilde L$,
one sees that these integrals depend continuously on $\Omega$ (for
$\Omega=0$, $E^*_{L,\cC,M_0}$ reduces to the energy $E_{L,M_0}$),
and the same is obviously true for the constraint manifold $M_0$. It is
known that, when $\Omega=0$, $J_1$, $J_2$ and $E^*_{L,\cC,M_0}$
are the components of a submersion from an open nonempty set
$M_\mathrm{reg}$ of $M_0$ to $\bR3$ \cite{zenkov,FGS2005}. Continuity
implies that, at least for $\Omega$ sufficiently close to zero, the
map $(J_1,J_2,E^*_{L,\cC,M_0})$ is a submersion from an open
nonempty subset of $M_\mathrm{reg}$ to $\bR3$.

2. When $\Omega=0$, the level curves of
$(J_1,J_2,E^*_{L,\cC,M_0})$  in $M_\mathrm{reg}$ are compact
\cite{zenkov,FGS2005}, hence bounded.  It follows that, for each
$\Omega$ sufficiently close to zero, there is an open nonempty subset
of $M_\mathrm{reg}$ where the submersion
$(J_1,J_2,E^*_{L,\cC,M_0})$ has bounded (hence compact, being a
submersion) level sets.

3. According to the mentioned reconstruction results from reduced periodic
dynamics by Field and Krupa, if the group is compact and
acts freely, then each `relative periodic orbit' (that is, the group
orbit in the phase space that projects over a periodic orbit of the
reduced system) is fibered by tori of dimension up to $1+\rho$, 
where $\rho$ is the rank of the group, on which motions are
quasi-periodic.
In our case, $\rho=2$.
\end{proof}

These are clearly partial results, that should be completed and extended
under several aspects. We shortly indicate some of them. 

First of all, the regions of the phase space where the dynamics is
quasi-periodic should be identified, and it should be understood
how they depend on $\Omega$ and on the shape of the surface.

A complementary question concerns the behaviour of motions that are not
quasi-periodic, if present. 

There are also interesting questions
about the geometry of the (singular) foliation by invariant tori.
The proof given above shows that each relative periodic orbit is
fibered by tori of some dimension between 1 and 3, but it does not
ensure that this dimension is the same across different relative
periodic orbits and that the
invariant tori are the fibers of a fibration of (an open subset of)
the phase space. This property is
important, because it implies the existence of the appropriate number
of first integrals that are usually associated to integrability (for
some results on this point in the case $\Omega=0$ see \cite{FG2007}).

Finally, even if the reconstruction procedure gives generically a fibration
by invariant tori of dimension 3, it might happen
that for certain shapes of the surface there are resonance conditions
among the frequencies of {\it all} 
the quasi-periodic motions---or equivalently, there exist
additional first integrals---and there is a fibration by
invariant tori of smaller dimension, either 2 (all motions have two
frequencies) or 1 (all motions are periodic). An instance of this
possibility is met in the limiting case of the sphere
on the turntable: since the $\mathrm{SO(3)}$-reduced system has
periodic dynamics, and $\mathrm{SO(3)}$ has rank~1, the unreduced
motions are quasi-periodic on tori of dimension at most~2.

This analysis (which has not yet been performed completely
for the case $\Omega=0$, either) requires manifestly 
an approach different from
the general one used in this section, and will be done elsewhere.

\begin{remark}{\rm 
The integrability result for the reduced system given 
in \cite{BMK2002}, based on Jacobi theorem,
could be strengthen if the common level sets of the
two first integrals $J_1$ and $J_2$ were compact. Under such a
hypothesis, Jacobi theorem implies that
these level sets are two-dimensional tori and motions
on them are linear
{\it after a time reparameterization}. 
However, not only the level sets of $J_1$ and $J_2$ are
unlikely to be compact (these two functions are linear in
some of the coordinates on the reduced phase space), but because of
the time reparameterization, this result
would be much weaker than those in Theorem~\ref{thm3}.
}\end{remark}

\vskip0.3cm
\noindent
{\bf Acknowledgements. } We thank Enrico Pagani for a useful
conversation and Larry Bates for suggesting the term `moving energy'.

\end{document}